\documentclass[12pt]{amsart}

\usepackage{setspace,a4wide,latexsym,amssymb,amsthm,amsmath}

\theoremstyle{plain}

\newtheorem{theorem}{Theorem}

\theoremstyle{remark}
\newtheorem{remark}{Remark}

\theoremstyle{definition}

\def\({\left (}
\def\){\right )}
\def\<{\left<}
\def\> { \right>}

\begin{document}

\title{Pseudo-parallel Lagrangian submanifolds are semi-parallel}

\author[F. Dillen]{Franki Dillen}
\author[J. Van der Veken]{Joeri Van der Veken}
\author[L. Vrancken]{Luc Vrancken}
\email[F. Dillen]{franki.dillen@wis.kuleuven.be} \email[J. Van der
Veken]{joeri.vanderveken@wis.kuleuven.be}\email[L.
Vrancken]{luc.vrancken@univ-valenciennes.fr}
\address[F. Dillen, J. Van der Veken, L. Vrancken]{Katholieke
Universiteit Leuven\\ Departement Wiskunde\\ Celestijnenlaan 200 B, Box 2400\\
BE-3001 Leuven\\ Belgium}
\address[L. Vrancken]{LAMAV \\ ISTV2\\ Universit\'e de Valenciennes\\
Campus du Mont Houy\\ 59313 Valenciennes Cedex 9\\ France}


\thanks{The second author is a postdoctoral researcher supported by the Research Foundation-Flanders (FWO)}
\thanks{Research supported by Research Foundation - Flanders project G.0432.07.}

\begin{abstract}  We prove a conjecture formulated by Pablo M. Chacon  and Guillermo A.
Lobos in \cite{CL} stating that every Lagrangian pseudo-parallel
submanifold of a complex space form of dimension at least 3 is
semi-parallel. We also propose to study another notion of
pseudo-parallelity which is more adapted to the Kaehlerian
setting.\end{abstract}

\keywords{Lagrangian submanifold, Pseudo-parallel submanifold,
Semi-parallel submanifold.}

\subjclass[2000]{Primary: 53C40.}

\maketitle

\section{Introduction.}
A submanifold $M^n$ is called pseudo-parallel if
\begin{equation}
\overline{R}\cdot h + \phi\, Q(g,h) = 0, \label{eqpseudopa}
\end{equation}
for some function $\phi$ on $M$, where is $\overline{R}$ is the
curvature operator of the Van der Waerden-Bortolotti connection, $h$
is the second fundamental form, $Q(g,h)$ is defined by
\begin{equation}\label{Q}
Q(g,h)(X,Y,U,V) = - ((X\wedge Y)\cdot h)(U,V) = h((X\wedge Y)U,V)
+h(U,(X\wedge Y)V)
\end{equation}
and $\wedge$ is defined by
\begin{equation*}
(X\wedge Y)Z  = g(Y,Z)X - g(X,Z)Y,
\end{equation*}
where $g$ is the metric. In formula (\ref{eqpseudopa}),
$\overline{R}$ acts on $h$ as a derivation, in particular
\begin{equation*}
(\overline{R}\cdot h)(X,Y,Z,W) =   R^{\perp}(X,Y)h(Z,W) -
h(R(X,Y)Z,W) - h(Z,R(X,Y)W),
\end{equation*}
where $R^{\perp}$ is the normal curvature tensor and $R$ is the
Riemann curvature tensor of $M$. Similarly in (\ref{Q}), $X\wedge Y$
acts as derivation on $h$.

Pseudo-parallel submanifolds are introduced in \cite{ALM1} and
\cite{ALM2} as generalization of semi-parallel submanifolds (i.e.
satisfying $\overline{R}\cdot h=0$), in the sense of \cite{DE}. The
notion pseudo-parallel generalizes semi-parallel in the same way as
pseudo-symmetry (in the sense of \cite{DES}) generalizes
semi-symmetry. See also \cite{DFHVV} for further discussions and
properties.

In this paper we consider pseudo-parallel Lagrangian submanifolds
$M^n$ of a complex space form $\widetilde M^n(4c)$ of constant
holomorphic curvature $4c$. In \cite{CL} it is proved that a
Lagrangian surface $M^2$ of a complex space form $\widetilde
M^2(4c)$ is pseudo-parallel if and only if it is flat or minimal and
the authors conjectured that for $n\geq 3$, pseudo-parallel implies
semi-parallel. We prove this conjecture.

\begin{theorem}
A Lagrangian pseudo-parallel submanifold $M^n$, $n\geq 3$, of a
complex space form $\widetilde M^n(4c)$ is semi-parallel.
\end{theorem}

Basically we only need the total symmetry of the cubic form $C$,
defined by $C(X,Y,Z) = g(h(X,Y),JZ)$ and the following identities
for the shape operator $A$ and the curvature tensors:
\begin{align*}
J(A_{JX}Y) &= h(X,Y) = J(A_{JY}X),\\
R^{\perp}(X,Y)JZ &= J R(X,Y)Z,\\
R(X,Y) &= c X\wedge Y + [A_{JX},A_{JY}].
\end{align*}

\section{Proof of the Theorem}

Let $M^n$ be pseudo-parallel. Suppose that $p\in M^n$ is a point
such that  $\phi(p) \ne 0$.

Since $g(h(V,W),JZ)$ is symmetric in $W$ and $Z$, from the Ricci
identity we immediately obtain that $g((\overline{R}\cdot
h)(X,Y,V,W),JZ)$ is symmetric in $W$ and $Z$. Since $M^n$ is
pseudo-parallel and $\phi(p) \ne 0$, we obtain that
$g(Q(g,h)(X,Y,V,W),JZ)$ is symmetric in $W$ and $Z$. This implies
that
\begin{equation*}
g(Y,W) g(h(V,X),JZ) - g(X,W) g(h(Y,V),JZ)
\end{equation*}
is symmetric in $W$ and $Z$, hence
$$
\begin{aligned}
g(Y,W) g(h(V,X),JZ) &- g(X,W) g(h(Y,V),JZ)\\ &= g(Y,Z) g(h(V,X),JW)
- g(X,Z) g(h(Y,V),JW).
\end{aligned}
$$
Taking $X=W=V$ and $Z,Y$ perpendicular to  $X$, we get
\begin{equation}\label{symmetry2}
- g(X,X) g(h(Y,Z),JX) =  g(Y,Z) g(h(X,X),JX).
\end{equation}

First we choose a unit vector $x$ of $T_pM$ such that the function
$f(v) = g(h(v,v),Jv)$, defined on the unit sphere of $T_pM$, attains
a maximal value at $x$. Then, by a standard argument, we obtain that
$g(h(x,x),Jv)=0$ for any unit vector $v\perp x$. Hence $A_{Jx}x =
\lambda_1 x$ for some number $\lambda_1$. Then taking $X=x$ in
(\ref{symmetry2}), we obtain that for any $y$ with $y\perp x$ holds
that $A_{Jx}y=-\lambda_1 y$.

For any such $y\perp x$, we take $X=y$ and $Y=Z=x$ in
(\ref{symmetry2}). This gives us
$$
g(h(y,y),Jy) = -g(y,y) g(h(x,x),Jy)=0,
$$
such that, by the total symmetry of the cubic form, $g(h(y,z),Jw)=0$
for all $y,z,w\perp x$. Hence $h$ has a very simple form:
\begin{equation}\label{form}
h(x,x)=\lambda_1 Jx,\quad h(x,y)=-\lambda_1 Jy, \quad h(y,z) =
-\lambda_1 g(y,z) Jx, \quad y,z\perp x.
\end{equation}
If $\lambda_1=0$, then $h$ vanishes at $p$. So we assume that
$\lambda_1\ne 0$.

Since $n>2$ we can choose orthonormal vectors $y,z$ orthogonal to
$x$. Then $R(x,y)y = (c-2\lambda_1^2)x$ and $R(y,z)z =
(c+\lambda_1^2)y$. Computing
$$
(\overline{R}\cdot h)(x,y,y,y) + \phi(p)\, Q(g,h)(x,y,y,y) = 0,
$$
gives us $3\lambda_1(c-2\lambda_1^2) = 2\phi(p) \lambda_1$ and hence
\begin{equation}\label{1}
\phi(p)= \frac{3}{2}(c-2\lambda_1^2).
\end{equation}
Similarly from
$$
(\overline{R}\cdot h)(x,y,y,z) + \phi(p)\, Q(g,h)(x,y,y,z) = 0,
$$
we obtain
\begin{equation}\label{2}
\phi(p)= (c-2\lambda_1^2).
\end{equation}
Hence from (\ref{1}) and (\ref{2}) we obtain that $\phi(p)=0$.

We conclude that either $\phi(p)=0$ or $h=0$ at $p$, in both cases
$\overline{R}\cdot h$ vanishes at $p$.\qed

\begin{remark}
In case $n=2$, then it follows immediately from (\ref{form}) that
pseudo-parallelity with $\phi(p)\ne 0$ implies minimality at $p$.
For $n>2$ equation (\ref{form}) implies that $M^n$ is $H$-umbilical
and we could have referred to \cite{CL} to obtain semi-parallelity,
but we finish the proof here for completeness.
\end{remark}

\begin{remark}
It seems that pseudo-parallelity is not the ideal condition to study
for Lagrangian submanifolds. In fact, the two sides of
(\ref{eqpseudopa}) don't have the same symmetries with respect to
$J$. It seems to be more interesting to study Lagrangian
submanifolds with pseudo-parallel cubic form, i.e. satisfying the
condition
\begin{equation}
R\cdot C + \phi\, Q(g,C) = 0.
\end{equation}
Since for a surface $R(X,Y)=K (X\wedge Y)$, where $K$ is the Gauss
curvature, it is clear that any Lagrangian surface has
pseudo-parallel cubic form with $\phi=K$.
\end{remark}

\end{document}